\def\draft#1{}                      

\documentclass[11pt]{article}
\usepackage{amssymb}

\def\eop{\hfill\framebox[2.4mm][t1]{\phantom{x}} \vskip 0.15cm } 

\title{On a theorem of Garza regarding \linebreak algebraic numbers with real conjugates}

\author{Gerald~H\"ohn\thanks{Department of Mathematics, 
Kansas State University, 138 Cardwell Hall, Manhattan, KS 66506-2602, USA.
E-mail: {\tt gerald@math.ksu.edu}}}

\begin{document}

\maketitle


\paragraph{1\ \ Introduction.}

For an algebraic number $\alpha$, that is, a root of an irreducible polynomial $\phi(x)$
with integer coefficients, the absolute
height of $\alpha$ is defined by $H(\alpha)=|c|^{1/d}\prod_{i=1}^d\max(1,\,|\alpha_i|)^{1/d}$
in case $\phi(x)=c\prod_{i=1}^d(x-\alpha_i)$.
The following lower estimate for the absolute height of~$\alpha$
was recently found by J.~Garza (\cite{Garza}, Theorem 1):

\medskip

\noindent{\bf Theorem:}\ \  {\em Let $\alpha\not=0$, $\pm 1$ be an algebraic number with $r>0$ 
real Galois conjugates. Then
$$ H(\alpha)\geq \left(\frac{2^{1-1/R}+\sqrt{4^{1-1/R}+4}}{2}\right)^{R/2}$$
where $R=r/d$ is the fraction of Galois conjugates $\alpha_i$ of $\alpha$ which are real.}

\medskip

If $R=1$, i.e., $\alpha$ is a totally real, the bound simplifies to Schinzel's
estimate (see~\cite{Schinzel}, Corollary~1')
$$ H(\alpha)\geq \left(\frac{1+\sqrt{5}}{2}\right)^{1/2}$$
stated in loc.~cit.~for algebraic integers only.
A short proof of Schinzel's bound in this case was given in~\cite{HoSchinzel}. In this
note we show that a similar method as in~\cite{HoSchinzel} together with basic properties
of absolute values of number fields also leads to a new derivation of Garza's bound.


\paragraph{2\ \ Proof of Theorem.} We start with an elementary estimate.

\smallskip 

\vbox{
\noindent{\bf Lemma:}\ \ {\em For $0<a<\frac{1}{2}$ let $f(x)=|x|^{1/2-a}|1-x^2|^a$. Then 
the function $f(x)/\max(1,\,|x|)$ has the global maximum $M_{\bf C}=2^a$ on the complex plane 
and the global maximum 
$$M_{\bf R}=(4a)^a(1-2a)^{1/4-a/2}(1+2a)^{-1/4-a/2}$$
on the real axis.}
}

{\it Proof of the lemma:} One has $f(x)\leq 2^a$ for $|x|\leq 1$ and $f(i)=2^a$. For $|x|\geq 1$
one gets $f(x)/|x|\leq |x|^{-1/2-a}(2|x|^2)^a\leq 2^a$ proving the
first statement. For the second statement, one verifies by using the first derivative and computing the
boundary values that $f(x)$ reaches the stated global maximum in the interval $[0,1]$ at 
$x_1=\sqrt{\frac{1-2a}{1+2a}}$ and that $f(x)/|x|$ reaches the same global maximum 
in the interval $[1,\infty)$ at $x_2=\sqrt{\frac{1+2a}{1-2a}}$. \eop

\smallskip

Continuing with the notation from the lemma, one has for an algebraic integer $\alpha$ the estimate
$$\prod_{i=1}^d f(\alpha_i)=|\phi(0)|^{1/2-a}|\phi(1)\phi(-1)|^a\geq 1.$$
Therefore,
$$\prod_{i=1}^d\max(1,\,|\alpha_i|)\geq  M_{\bf R}^{-r} M_{\bf C}^{r-d} \prod_{i=1}^d  f(\alpha_i) \geq  M_{\bf R}^{-r} M_{\bf C}^{r-d}$$
or $H(\alpha)\geq  M_{\bf R}^{-R} M_{\bf C}^{R-1}$ for the height. 
Applying the lemma for $a=\frac{1}{2}(1+4^{1/R})^{-1/2}$ gives
\begin{eqnarray*}
H(\alpha) & \geq  & (4a)^{-aR}(1-2a)^{(a/2-1/4)R}(1+2a)^{(a/2+1/4)R}2^{a(R-1)} \\
 &= & \left( \left(\frac{1+4^{1/R}}{4}\right)^a\left(\frac{4^{1/R}}{1+4^{1/R}}\right)^a4^{a(1-1/R)}\cdot  \frac{1+2a}{(1-4a^2)^{1/2}} \right)^{R/2}    \\
 &= & \left( \left(\frac{1+4^{1/R}}{4^{1/R}}\right)^{1/2}\left(1+(1+4^{1/R})^{-1/2}\right)  \right)^{R/2} =
 \left(\frac{2^{1-1/R}+\sqrt{4^{1-1/R}+4}}{2}\right)^{R/2}, 
\end{eqnarray*}
which finishes the proof of the theorem in the case of the algebraic integers.

\smallskip

The above argument can be extended to abitrary algebraic numbers $\alpha$
by using some basic algebraic number theory and properties of the absolute height 
(cf.~\cite{Ishak} for the case of Schinzel's result).

Let $k={\bf Q}(\alpha)$. For a place $\nu$ of $k$ we denote by  $|\,.\,|_{\nu}$ the corresponding normalized 
absolute value of $k$, so that $\prod_{v}|\beta|_\nu=1$ for a non-zero algebraic number $\beta$ in $k$. 
Then the absolute height of $\beta$  equals  
$H(\beta)=\prod_{v}\max(1,|\beta|_\nu)$. With $a\leq 1/2$ as above, we have the estimate
\begin{eqnarray*}
1 & = & \prod_{\nu}|\alpha-\alpha^{-1}|_\nu^a = \prod_{\nu|\infty}|\alpha-\alpha^{-1}|_\nu^a \cdot 
                          \prod_{\nu\nmid\infty}| \alpha-\alpha^{-1} |_\nu^a \\
  & \leq & \prod_{\nu|\infty} |\alpha-\alpha^{-1}|_\nu^a
\prod_{\nu\nmid\infty}  \max(1,|\alpha|_\nu)^{a} \max(1,|\alpha^{-1}|_\nu)^{a}\\
  & \leq    & \prod_{\nu|\infty} \frac{\big(|\alpha_\nu-{\alpha_\nu}^{-1}|^a\big)^{d_\nu/d}}
                       {\big(\max(1,|\alpha_\nu|)^{1/2}\max(1,|\alpha_\nu^{-1}|)^{1/2}\big)^{d_\nu/d}} \cdot 
                                     \prod_{\nu} \max(1,|\alpha|_\nu)^{1/2} \max(1,|\alpha|^{-1}_\nu)^{1/2} 
\end{eqnarray*}
where $d_\nu=[\,k_\nu\,{:}\,{\bf R}]$ and $\alpha_\nu$ is the image of $\alpha$ under some Galois automorphism of the Galois closure of $k$ 
such that $|\alpha|_\nu=|\alpha_\nu|^{d_\nu/d}=|\alpha_i|^{d_\nu/d}$ for some $i$ so that one factor for each pair 
$\{\alpha_i,\,\bar\alpha_i\}$ appears in the product over the 
archimedean places.
Since $g(x)=|x-x^{-1}|^a/(\max(1,|x|)^{1/2}\max(1,|x^{-1}|)^{1/2})$ is symmetric under $x\mapsto x^{-1}$ we can assume $|x|\geq 1$ where
$g(x)=f(x)/\max(1,|x|)$. By applying the lemma we get now the estimate
$$1\leq  M_{\bf R}^{R} M_{\bf C}^{1-R} \cdot H(\alpha)^{1/2}H(\alpha^{-1})^{1/2}$$
and the result follows as before by using $H(\alpha)=H(\alpha^{-1})$.


\paragraph{3\ \ Remarks.} 1. Under all functions $\tilde f(x)=|x|^u|1-x^2|^v$, the chosen $f(x)$ gives
the best estimate for $H(\alpha)$.

\smallskip

2. For $R=1$ the bound for $H(\alpha)$ is optimal. One may ask if this is also the case for other values of $R$, although
it follows from the proof that there cannot exist an $\alpha$ actually reaching the bound.

\smallskip

3. The main difference to Garza's proof is that we replace 
a sequence of inequalities in~\cite{Garza} with the estimate of the lemma,
allowing a particular elementary proof for algebraic integers.


\end{document}